\documentclass[a4paper]{article}
\usepackage{amssymb}
\usepackage{amsmath}

\providecommand{\dsum}{\displaystyle\sum}
\providecommand{\TEXTsymbol}[1]{}

\begin{document}

\title{A Classification of Self-Maps of Generalized Grassmannians}
\author{Haibao Duan, Ruizhi Huang, Xuezhi Zhao}
\date{ }
\maketitle

\begin{abstract}
Generalized Grassmannians form a fundamental class of flag manifolds
associated with Lie groups. The purpose of this paper is to classify
self-maps of generalized Grassmannians of nonzero degree in terms of their
induced actions on cohomology.

We prove a rigidity theorem showing that, despite the rich and intricate
structure of their cohomology rings, the induced cohomology endomorphisms
fall into only two natural types: Adams operations determined by the degree
and Dynkin symmetries arising from automorphisms of the Dynkin diagram.

Combining the geometry of root and weight systems, the actions of Weyl and
Dynkin symmetries on Schubert classes, and Bott--Samelson desingularizations
of Schubert varieties, our approach applies uniformly to generalized
Grassmannians of all Lie types.
\end{abstract}

\section{Introduction}

Let $G$ be a compact connected semisimple Lie group with Lie algebra $L(G)$
and exponential map $\exp :L(G)\rightarrow G$. For a nonzero element $\alpha
\in L(G)\backslash \left\{ 0\right\} $ let $P_{\alpha }$ denote the
centralizer of the one-parameter subgroup

\begin{center}
$\left\{ \exp (t\alpha )\in G\mid t\in \mathbb{R}\right\} $.
\end{center}

\noindent Then $P_{\alpha }$ is a parabolic subgroup, and the homogeneous
space $G/P_{\alpha }$ is canonically a smooth complex projective variety,
called a \textsl{flag manifold} of $G$. These manifolds form a fundamental
class of compact homogeneous spaces. Their cohomology has played a central
role in the development of characteristic classes, representation theory,
and algebraic geometry, and has been extensively studied through Schubert
calculus; see, for example, \cite{BGP,DZ4}. 

On the other hand, for a topological space $X$, let $[X,X]$ denote the
monoid of homotopy classes of self-maps of $X$, and let $\mathrm{End}%
(H^{\ast }(X))$ denote the monoid of endomorphisms of the integral
cohomology ring $H^{\ast }(X)$. Every self-map $f$ of $X$ induces a graded
ring endomorphism $f^{\ast }$, giving rise to the natural representation

\begin{enumerate}
\item[(1.1)] $R:[X,X]\rightarrow \mathrm{End}(H^{\ast }(X))$, $R[f]=f^{\ast
} $.
\end{enumerate}

\noindent When $X=G/P$ is a flag manifold, this representation is known to
be faithful up to finite ambiguity \cite{HM,Sh,P}. More precisely, let $X_{0}
$ denote the rationalization of $X$. By combining Sullivan's rational
homotopy theory \cite{S} with the formality of flag manifolds, one obtains
the following result; see, for example, \cite{GH,P}.

\bigskip

\noindent \textbf{Lemma 1.1. }If $X$ is a flag manifold, the rationalization
of $R$ induces a bijection

\begin{center}
$[X_{0},X_{0}]\rightarrow \mathrm{End}(H^{\ast }(X;\mathbb{Q}))$.
\end{center}

Lemma 1.1 reduces the topological problem of determining the homotopy set $%
[X,X]$ of a flag manifold largely to the algebraic problem of describing the
monoid $\mathrm{End}(H^{\ast }(X))$, and has motivated a substantial body of
work in this direction; see, for examples, \cite{D2,DZ,GH,H1,KT,Lin,P} or
Remark 1.7. The aim of this paper is to make the representation $R$ explicit
for an important family of flag manifolds, namely, the generalized
Grassmannians, which we now define.

In terms of the root system of $G$, flag manifolds admit the following
classification \cite[Lemma 2.2]{DZ2}. Fix a maximal torus $T\subset G$, and
let

\begin{center}
$\Omega =\{\omega _{{1}},\ldots ,\omega _{{n}}\}\subset L(T)$
\end{center}

\noindent be a set of \textsl{fundamental dominant weights}; see \cite[p.67]%
{Hu} or Section 3.1. For each subset $I\subseteq \{1,\cdots ,n\}$, let $P_{I}
$ denote the centralizer of the one-parameter subgroup

\begin{center}
$\alpha :\mathbb{R}\rightarrow G$, $\alpha (t)=\exp (t\sum\limits_{i\in
I}\omega _{i})$.
\end{center}

\noindent \textbf{Lemma 1.2}.\textbf{\ }Every parabolic subgroup $P$ of $G$
is conjugate to $P_{I}$\ for some $I\subseteq \{1,\cdots ,n\}$.
Consequently, every flag manifold of $G$ is isomorphic to a homogeneous
space $G/P_{I}$.

\bigskip

Among the flag manifolds described in Lemma 1.2, two families will be of
particular importance.

(a) If $I=\{1,\cdots ,n\}$, then $P_{I}=T$, and $G/P_{I}=G/T$ is the\textsl{%
\ complete flag manifold} of $G$.

(b) If $I=\{i\}$, then $G/P_{\{i\}}$ is called the \textsl{generalized} 
\textsl{Grassmannian\ }associated with the fundamental weight $\omega _{i}$.

Generalized Grassmannians satisfy the following basic properties, which will
be used repeatedly throughout the paper; see \cite[Lemma 2.2]{DZ2}.

\bigskip

\noindent \textbf{Lemma 1.3}.\textbf{\ }Let $X$ be a generalized
Grassmannian. Then

(i) $X$ is isomorphic to $G_{0}/P_{\{i\}}$, where $G_{0}$ is a simply
connected simple Lie group;

(ii) the integral cohomology $H^{\ast }(X)$ is torsion free and concentrated
in the even degrees;

(iii) $H^{2}(X)\cong \mathbb{Z}$, generated by the canonical K\"{a}hler
class $\kappa _{X}$ of $X$.

\bigskip

In view of Lemma 1.3(i), we can assume, without loss of generality, that the
ambient Lie group $G$ is simply connected and simple. Such Lie groups are
classified by their Dynkin diagrams $\Gamma _{G}$, which consists of the
four infinite classical families

\begin{center}
$\left\{ A_{l},\ l\geq 1\right\} $, $\left\{ B_{l},\ l\geq 2\right\} $, $%
\left\{ C_{l},\ l\geq 3\right\} ,\left\{ D_{l},\ l\geq 4\right\} $,
\end{center}

\noindent together with the five exceptional types $G_{2},F_{4},E_{6},E_{7}$%
, $E_{8}$ (see \cite[p.58]{Hu}). Let $\mathrm{Aut}(\Gamma _{G})$ denote the
automorphism group of the Dynkin diagram $\Gamma _{G}$. Throughout, we index
the fundamental dominant weights in $%
\Omega
$ according to the vertices of $\Gamma _{G}$, following the convention in 
\cite[p.58]{Hu}. Since $\mathrm{Aut}(\Gamma _{G})$ acts on $\Omega $ by
permuting the fundamental dominant weights, for each $\omega _{{i}}\in
\Omega $ we define its stabilizer by

\begin{center}
$\mathrm{Aut}(\Gamma _{G},\omega _{{i}}):=\left\{ \sigma \in \mathrm{Aut}%
(\Gamma _{G})\mid \sigma (\omega _{{i}})=\omega _{{i}}\right\} $.
\end{center}

\noindent A direct inspection of the Dynkin diagrams gives the following
classification.

\bigskip

\noindent \textbf{Lemma 1.4.} The stabilizer $\mathrm{Aut}(\Gamma
_{G},\omega _{{i}})$ is nontrivial if and only if the pair $(G,i)$ appears
in the following table.

\begin{center}
\begin{tabular}{l||l|l|l|l}
\hline
$G$ & $A_{2k+1}$ & $D_{l}$, $l>4$ & $D_{4}$ & $E_{6}$ \\ \hline
$i$ & $k+1$ & $1,\cdots ,n-2$ & $2,1,3,4$ & $2,4$ \\ \hline
$\mathrm{Aut}(\Gamma _{G},\omega _{{i}})$ & $\mathbb{Z}_{2}$ & $\mathbb{Z}%
_{2}$ & $S_{3},S_{2},S_{2},S_{2}$ & $\mathbb{Z}_{2}$ \\ \hline
\end{tabular}%
.
\end{center}

\noindent Here $S_{m}$ denotes the symmetric group on $m$ letters.

\bigskip

The nontrivial stabilizers classified in Lemma 1.4 give rise to nontrivial
diffeomorphisms of generalized Grassmannians and, consequently, to
automorphisms of their cohomology rings. More precisely, Lemma 3.2 shows
that every automorphism $\sigma \in \mathrm{Aut}(\Gamma _{G},\omega _{{i}})$
induces a diffeomorphism

\begin{center}
$h_{\sigma }:G/P_{\{i\}}\rightarrow G/P_{\{i\}}$
\end{center}

\noindent and subsequently, gives rise to a homomorphism

\begin{center}
$T:\mathrm{Aut}(\Gamma _{G},\omega _{{i}})\rightarrow \mathrm{Aut}(H^{\ast
}(G/P_{\{i\}}))$,$\quad T(\sigma )=h_{\sigma }^{\ast }$.
\end{center}

\noindent Our first main result ensures that this representation is always
faithful.

\bigskip

\noindent \textbf{Theorem A. }For every generalized Grassmannian $%
G/P_{\{i\}} $, the map $T$ is an injective homomorphism. In particular,

\begin{center}
$\sigma \neq \mathrm{id}\Longrightarrow h_{\sigma }^{\ast }\neq \mathrm{Id}$.
\end{center}

Let $X=G/P_{\{i\}}$ be a generalized Grassmannian, and let $\sigma \in 
\mathrm{Aut}(\Gamma _{G},\omega _{{i}})$ be nontrivial. The automorphism $%
h_{\sigma }^{\ast }\in \mathrm{Aut}(H^{\ast }(X))$ arising in Theorem A is
called a \textsl{Dynkin symmetry} of the ring $H^{\ast }(X)$.

Since $H^{\ast }(X)$ is torsion free and concentrated in even degrees, for
any nonzero integer $k$, the \textsl{Adams operator }$\psi ^{k}\in \mathrm{%
End}(H^{\ast }(X))$ is well defined by

\begin{center}
$\psi ^{k}(x)=k^{r}x,$ $\quad x\in H^{2r}(X)$.
\end{center}

For a self-map $f:X\rightarrow X$, since $H^{2}(X)\cong \mathbb{Z}$, we
define the degree $\deg (f)\in \mathbb{Z}$ by

\begin{center}
$f^{\ast }(\kappa _{X})=\deg (f)\cdot \kappa _{X}$.
\end{center}

\noindent Our second main theorem establishes a rigidity phenomenon:
whenever $\deg (f)\neq 0$, the induced action $f^{\ast }$ on the integral
cohomology ring is uniquely determined by the integer $\deg (f)$, up to
composition with a Dynkin symmetry.

\bigskip

\noindent \textbf{Theorem B. }Let $f$ be a self-map of a generalized
Grassmannian $X$ with $\deg (f)\neq 0$. Then

\begin{center}
$f^{\ast }=\psi ^{\deg (f)}$ or $f^{\ast }=\psi ^{\deg (f)}\circ \tau $,
\end{center}

\noindent where the second possibility can occur only if $\mathrm{Aut}%
(\Gamma _{G},\omega _{{i}})\neq \{\mathrm{id}\}$, and $\tau $ is a Dynkin
symmetry of the ring $H^{\ast }(X)$ associated with a nontrivial element $%
\sigma \in \mathrm{Aut}(\Gamma _{G},\omega _{{i}})$ as classified in Lemma
1.4.

\bigskip

Nontrivial Dynkin symmetries occur only in exceptional cases among
generalized Grassmannians. For example, when $G$ is an exceptional Lie
group, Lemma 1.4 shows that $\mathrm{Aut}(\Gamma _{G},\omega _{{i}})=\left\{ 
\mathrm{id}\right\} $ with the only exceptions

\begin{center}
$\mathrm{Aut}(\Gamma _{E_{6}},\omega _{{2}})\cong \mathrm{Aut}(\Gamma
_{E_{6}},\omega _{{4}})\cong \mathbb{Z}_{2}$.
\end{center}

\noindent This immediately gives the following consequence of Theorem B.

\bigskip 

\noindent \textbf{Corollary 1.5. }Let $X=G/P_{\{i\}}$ be a generalized
Grassmannian with $G$ of exceptional Lie type, and suppose that 

\begin{center}
$(G,i)\neq (E_{6},2)$, $(E_{6},4)$. 
\end{center}

\noindent Then for every self-map $f$ of $X$ with $\deg (f)\neq 0$,

\begin{center}
$f^{\ast }=\psi ^{\deg (f)}$.
\end{center}

For a simply connected finite $CW$ complex $X$, let $\mathrm{Aut}(X)$ denote
the group of self-homotopy equivalences of $X$. If $X$ is a generalized
Grassmannian and $f\in \mathrm{Aut}(X)$, then necessarily

\begin{center}
$\deg (f)=\pm 1$.
\end{center}

\noindent Theorem B therefore gives the following immediate classification
of the actions of self-homotopy equivalences on the cohomology of a
generalized Grassmannian.

\bigskip

\noindent \textbf{Corollary 1.6. }For every $f\in \mathrm{Aut}(X)$,

\begin{center}
$f^{\ast }=\psi ^{\pm 1}$ or $f^{\ast }=\psi ^{\pm 1}\circ \tau $,
\end{center}

\noindent where $\tau $ is a Dynkin symmetry of the ring $H^{\ast }(X)$.

\bigskip

\noindent \textbf{Remark 1.7.} Motivated by Lemma 1.1, Hoffman \cite{H1}
established Theorem B for the complex Grassmannian $X=G_{n,r}(\mathbb{C})$
of $r$-dimensional complex subspaces of $\mathbb{C}^{n}$, while the first
author \cite{D2} proved the corresponding result for the Grassmannian $X=%
\mathbb{C}S_{n}$ of complex structures on the Euclidean space $\mathbb{R}%
^{2n}$. Both arguments rely on explicit presentations of the corresponding
cohomology rings. However, this approach does not readily extend to the
general setting, as it would inevitably involve case-by-case computations:
by Lemma 1.2, a simple Lie group of rank $n$ gives rise to $n$ generalized
Grassmannians, whose cohomology rings may vary substantially with the choice
of the parabolic subgroup $P_{\{i\}},1\leq i\leq n$; see, for example, \cite[%
Theorems 1-7]{DZ1}.

In contrast, the proof of Theorem B exploits structural properties common to
all generalized Grassmannians, yielding a uniform argument across all Lie
types without requiring explicit presentations of their individual
cohomology rings.

\bigskip

The remainder of the paper is organized as follows. In Section 2, using the
hard Lefschetz theorem for complex projective manifolds, we show that every
self-map of a generalized Grassmannian with nonzero degree is a rational
homotopy equivalence. In Section 3, we recall two types of symmetries of the
complete flag manifold $G/T$ and determine which of them descend to
nontrivial symmetries of the generalized Grassmannian $G/P_{\{i\}}$. Section
4 develops the necessary consequences of the basis theorem in Schubert
calculus. Combining these ingredients, we complete the proofs of Theorems A
and B in Section 5.

For a generalized Grassmannian $X$ the set $[X,X]$ admits the partition

\begin{center}
$[X,X]=\amalg _{k\in \mathbb{Z}}[X,X]_{k}$,
\end{center}

\noindent where

\begin{center}
$[X,X]_{k}:=\left\{ [f]\in \lbrack X,X]\mid \deg (f)=k\right\} $.
\end{center}

\noindent Theorem B naturally leads to two further questions:

\begin{enumerate}
\item[(1)] Determine the degree-zero component $[X,X]_{0}$ of $[X,X]$;

\item[(2)] For which nonzero integer $k$, is the Adams operator\textsl{\ }$%
\psi ^{k}$ realizable by a self-map $f$ of $X$?
\end{enumerate}

\noindent As complements to Theorem B, we provide partial answers to these
questions in Section 6.

Throughout the paper, the cohomology of a topological space is understood to
have coefficients in $\mathbb{Z}$, unless otherwise specified.

\section{Applications of the hard Lefschetz theorem to homotopy theory}

Throughout this section, let $M$ be a simply connected closed complex
projective manifold of complex dimension $n$, and let $\omega _{M}\in
H^{2}(M)$ denote an integral K\"{a}hler class. The hard Lefschetz theorem 
\cite[p.122]{GH1} asserts that, for every $0\leq r\leq n$, the homomorphism

\begin{enumerate}
\item[(2.1)] $L^{n-r}:H^{r}(M;\mathbb{R})\rightarrow H^{2n-r}(M;\mathbb{R})$%
, $x\rightarrow \omega _{M}^{n-r}\cup x$,
\end{enumerate}

\noindent is an isomorphism. We say that a self-map $f:M\rightarrow M$ 
\textsl{preserves the K\"{a}hler class} if

\begin{center}
$f^{\ast }(\omega _{M})=c\cdot \omega _{M}$
\end{center}

\noindent for some nonzero integer $c$. The main result of this section is
the following.

\bigskip

\noindent \textbf{Theorem 2.1}. Let $f$ be a self-map of $M$ that preserves
the K\"{a}hler class. Then $f$ is a rational homotopy equivalence.

\bigskip

\noindent \textbf{Proof. }The proof proceeds in two steps.

\textbf{Step 1. }For every $0\leq r\leq n$, the induced homomorphism $%
f^{\ast }:H^{r}(M;\mathbb{R})\rightarrow H^{r}(M;\mathbb{R})$ is an
isomorphism.

Fix an integer $r$ with $0\leq r\leq n$, let $b_{r}=\dim H^{r}(M;\mathbb{R})$%
, and choose a basis

\begin{center}
$\{x_{1},\cdots ,x_{b_{r}}\}$
\end{center}

\noindent of $H^{r}(M;\mathbb{R})$. Since $H^{2n}(M;\mathbb{R})\cong \mathbb{%
R}$ is generated by $\omega _{M}^{n}$, there exists a unique matrix $%
A_{r}=(a_{ij})\in M_{b_{r}}(\mathbb{R})$ such that

\begin{enumerate}
\item[(2.2)] \noindent $\omega _{M}^{n-r}\cup x_{i}\cup x_{j}=a_{ij}\omega
_{M}^{n}$, \quad $1\leq i,j\leq b_{r}$.
\end{enumerate}

\noindent By (2.1), the set

\begin{center}
$\{\omega _{M}^{n-r}\cup x_{i},1\leq i\leq b_{r}\}$
\end{center}

\noindent is a basis of $H^{2n-r}(M;\mathbb{R})$. Thus $A_{r}$ represents
the Poincar\'{e} duality pairing

\begin{center}
$H^{r}(M;\mathbb{R})\times H^{2n-r}(M;\mathbb{R})\rightarrow H^{2n}(M;%
\mathbb{R})\cong \mathbb{R}$
\end{center}

\noindent and hence is nonsingular.

Assume that $f^{\ast }(\omega _{M})=c\cdot \omega _{M}$, $c\neq 0$. Writing

\begin{enumerate}
\item[(2.3)] $f^{\ast }(x_{i})=\dsum\limits_{j=1}^{b_{r}}h_{ij}x_{j}$, $%
B_{r}(f)=(h_{ij})\in M_{b_{r}}(\mathbb{R})$,
\end{enumerate}

\noindent it suffices to prove that $B_{r}(f)$ is non-singular.

Applying $f^{\ast }$ to (2.2) and using the naturality of the cup product
together with (2.3), we obtain

\begin{center}
$c^{n-r}B_{r}(f)^{T}A_{r}B_{r}(f)=c^{n}A_{r}$,
\end{center}

\noindent where $B_{r}(f)^{T}$ denotes the transpose of $B_{r}(f)$. Taking
determinants and using the nonsingularity of $A_{r}$, we obtain

\begin{enumerate}
\item[(2.4)] $\det (B_{r}(f))^{2}=c^{r\cdot b_{r}}\neq 0$.
\end{enumerate}

\noindent Therefore $B_{r}(f)$ is invertible, proving Step 1.

\bigskip

\textbf{Step 2.} The map $f$ is a rational homotopy equivalence.

By Step 1, $f^{\ast }$ is an isomorphism for every $r\leq $ $n$. Poincar\'{e}
duality then implies that $f^{\ast }$ is an isomorphism in every degree.
Hence

\begin{center}
$f^{\ast }:$ $H^{\ast }(M;\mathbb{Q})\rightarrow H^{\ast }(M;\mathbb{Q})$
\end{center}

\noindent is a graded ring isomorphism.

Since $M$ has the homotopy type of a simply connected finite $CW$ complex,
the rational Whitehead theorem (see Sullivan \cite{S} or Griffiths--Morgan 
\cite[Ch.5]{GM}) implies that $f$ is a rational homotopy equivalence.\hfill $%
\square $

\bigskip

By Lemma 1.3(iii), every self-map of a generalized Grassmannian preserves
its K\"{a}hler class. The following consequence of Theorem 2.1 will play a
fundamental role in what follows.

\bigskip

\noindent \textbf{Corollary 2.2. }Let $f:X\rightarrow X$ be a self-map of a
generalized Grassmannian $X=G/P_{\{i\}}$ with $\deg (f)\neq 0$. Then $f$ is
a rational homotopy equivalence.\hfill $\square $

\bigskip

Let $X$ be a simply connected finite $CW$ complex. The classical Whitehead
theorem implies that the group of self-homotopy equivalences of $X$ is given
by

\begin{enumerate}
\item[(2.5)] $\mathrm{Aut}(X)=\left\{ [f]\in \lbrack X,X]\mid f^{\ast }\in 
\mathrm{Aut}(H^{\ast }(X))\right\} $.
\end{enumerate}

\noindent In contrast, the following integral refinement of Corollary 2.2
provides a simple characterization of the self-homotopy equivalences of
generalized Grassmannians.

\bigskip

\noindent \textbf{Theorem 2.3. }Let $X$ be a generalized Grassmannian. Then

\begin{center}
$\mathrm{Aut}(X)=\left\{ [f]\in \lbrack X,X]\mid \deg (f)=\pm 1\right\} $.
\end{center}

\noindent \textbf{Proof. }Since $H^{2}(X)\cong \mathbb{Z}$ is generated by
the K\"{a}hler class $\kappa _{X}$ (Lemma 1.3), it suffices to prove that
every self-map $f:X\rightarrow X$ with $\deg (f)=\pm 1$ is a homotopy
equivalence.

Since $H^{\ast }(X)$ is torsion free (Lemma 1.3), the basis $\{x_{1},\cdots
,x_{b_{r}}\}$ in the proof of Theorem 2.1 may be chosen in integral
cohomology. Thus,

\begin{center}
$B_{r}(f)=(h_{ij})\in M_{b_{r}}(\mathbb{Z})$.
\end{center}

\noindent Since $\deg (f)=\pm 1$, (2.4) gives

\begin{center}
$\det B_{r}(f)^{2}=1$.
\end{center}

\noindent Hence $B_{r}(f)$ is unimodular, and $f^{\ast }$ is an automorphism
of $H^{r}(X)$ for every $r\leq n$. Poincar\'{e} duality then implies that $%
f^{\ast }$ is an automorphism of $H^{\ast }(X)$. The conclusion follows from
(2.5).\hfill $\square $

\section{The Weyl and Dynkin symmetries on $G/T$}

Throughout this section, let $G$ be a simply connected simple Lie group with
maximal torus $T$, and let $P_{\left\{ i\right\} }$ be the parabolic
subgroup associated with the fundamental weight $\omega _{i}$. The inclusions

\begin{center}
$T\subset P_{\{i\}}\subset G$
\end{center}

\noindent give rise to the canonical fibration

\begin{enumerate}
\item[(3.1)] $P_{\{i\}}/T\hookrightarrow G/T\overset{\pi }{\rightarrow }%
G/P_{\{i\}}$,
\end{enumerate}

\noindent relating the complete flag manifold $G/T$ to the generalized
Grassmannian $G/P_{\{i\}}$.

The complete flag manifold $G/T$ admits two natural types of symmetries
associated with the root system of $G$: Weyl symmetries arising from the
Weyl group $W_{G}$, and Dynkin symmetries arising from automorphisms of the
Dynkin diagram $\Gamma _{G}$. In this section we determine which of these
symmetries descend, via the projection $\pi $, to nontrivial symmetries of
the generalized Grassmannian $G/P_{\{i\}}$.

Theorem 3.3 recalls a result of Papadima, describing the automorphism group $%
\mathrm{Aut}(H^{\ast }(G/T))$ in terms of these two types of geometric
symmetries of $G/T$. The main result in this section is Theorem 3.4, which
provides the missing link between Theorem 3.3 and the proof of Theorem B.

We begin with a geometric introduction to the root and weight systems of $G$%
, thereby fixing the notation used throughout the remainder of the paper and
developing a geometric framework for the arguments in the sequel.

\subsection{A geometric introduction to the root and weight systems}

The Lie algebra $L(G)$ (resp. the Cartan subalgebra $L(T)$) is identified
with the tangent space to $G$ (resp. $T$) at the identity $e$. Equip the
real vector space $L(G)$ with a fixed \textrm{Ad}$_{G}$-invariant inner
product $(,)$, so that

(a) the exponential map of $G$ at the unit $e$ is well defined, and yields
the commutative diagram

\begin{center}
\begin{tabular}{lll}
$\quad L(T)$ & $\rightarrow $ & $L(G)$ \\ 
$\exp \downarrow $ &  & $\downarrow \exp $ \\ 
$\quad T$ & $\rightarrow $ & $G$%
\end{tabular}%
,
\end{center}

\noindent where the horizontal arrows are the natural inclusions.

(b) the invariant inner product identifies $L(T)$ with its dual space $%
L(T)^{\ast }$.

Set $m:=\frac{1}{2}(\dim G-\dim T)$. The singular set of $\exp $ determines
a finite collection of $m$ hyperplanes through the origin in $L(T)$, called
the singular hyperplanes of $G$ \cite[p.168]{BT}:

\begin{center}
$\mathcal{S}(G)=\{L_{1},\cdots ,L_{m}\}$.
\end{center}

\noindent Moreover, in view of the geometric fact that the exponential map
carries each normal line $l_{k}$ to $L_{k}$ through the origin $0\in L(T)$
onto a circle subgroup on $T$, for each singular hyperplane $L_{k}$, let $%
\pm \alpha _{k}\in l_{k}$ be the pair of shortest nonzero vectors satisfying

\begin{center}
$\exp (\pm \alpha _{k})=e$, $1\leq k\leq m$.
\end{center}

\noindent These vectors gives rise to the \textsl{root system of }$G$

\begin{center}
$\Phi :=\{\pm \alpha _{k}\in L(T)\mid 1\leq k\leq m\}$.
\end{center}

In addition, the planes in $\mathcal{S}(G)$ divide $L(T)$ into finitely many
convex regions, each of them is called a \textsl{Weyl chamber} of $G$. Fix a
regular point $x_{0}\in L(T)$, and denote by $\mathcal{F}(x_{0})$ the
closure of the Weyl chamber containing $x_{0}$. Let

\begin{enumerate}
\item[(3.2)] $L(x_{0})=\{L_{1},\cdots ,L_{n}\}\subset \mathcal{S}(G)$, $%
n=\dim T$,
\end{enumerate}

\noindent denote the collection of walls of $\mathcal{F}(x_{0})$, and let $%
\alpha _{i}\in \Phi $ be the unique root orthogonal to $L_{i}$ and pointing
toward the interior of $\mathcal{F}(x_{0})$. The subset

\begin{center}
$\Delta =\{\alpha _{1},\cdots ,\alpha _{n}\}\subset \Phi $
\end{center}

\noindent so obtained is the set of\textsl{\ simple roots} of $G$ relative
to the chamber $\mathcal{F}(x_{0})$ \cite[p.49]{Hu}. The \textsl{fundamental
dominant weights} relative to $\mathcal{F}(x_{0})$ are the vectors

\begin{center}
$\Omega =\{\omega _{1},\cdots ,\omega _{n}\}\subset L(T)$
\end{center}

\noindent characterized by the linear system \cite[p.49]{Hu}

\begin{center}
$2(\omega _{i},\alpha _{j})/(\alpha _{j},\alpha _{j})=\delta _{ij},$ $1\leq
i,j\leq n$.
\end{center}

\noindent Equivalently \cite[Lemma 2.3(iii)]{D3},

\begin{center}
$A\left( 
\begin{tabular}{l}
$\omega _{1}$ \\ 
$\vdots $ \\ 
$\omega _{n}$%
\end{tabular}%
\right) =\left( 
\begin{tabular}{l}
$\alpha _{1}$ \\ 
$\vdots $ \\ 
$\alpha _{n}$%
\end{tabular}%
\right) $,
\end{center}

\noindent where $A\in M_{n}(\mathbb{Z})$ is the \textsl{Cartan matrix} of $G$
\cite[p.55]{Hu}. Over the ring of integers, the set $\Omega =\{\omega
_{1},\cdots ,\omega _{n}\}$ generates a lattice $\Lambda $, called \textsl{%
the weight lattice} of $G$.

Geometrically, the fundamental weight $\omega _{i}$ spans the unique edge

\begin{center}
$\left\{ t\omega _{i}\in L(T)\mid t\geq 0\right\} $
\end{center}

\noindent of the chamber $\mathcal{F}(x_{0})$ opposite to the wall $L_{i}$.
This geometric interpretation explains why the fundamental dominant weights
play a distinguished role in the classification of parabolic subgroups and
generalized Grassmannians (see Lemma 1.2).

\subsection{The Weyl symmetries on $G/T$}

Recall that the \textsl{Weyl group} of $G$ (with respect to a fixed maximal
torus $T$) is defined by

\begin{center}
$W_{G}:=N_{G}(T)/T$,
\end{center}

\noindent where $N_{G}(T)$ denotes the normalizer of $T$

\begin{center}
$N_{G}(T)=\{g\in G\mid gTg^{-1}=T\}$.
\end{center}

\noindent The right action of $N_{G}(T)$ on $G/T$ induces an injective
homomorphism

\begin{enumerate}
\item[(3.3)] $W_{G}=N_{G}(T)/T\hookrightarrow \mathrm{Diff}(G/T)$.
\end{enumerate}

\noindent The diffeomorphisms in the image of (3.3) are called \textsl{Weyl
symmetries} of $G/T$.

Under the canonical action of $W_{G}$ on the Cartan subalgebra $L(T)$, the
Weyl group $W_{P_{\{i\}}}$ of $P_{\{i\}}$ is naturally identified with the
stabilizer of $\omega _{i}$

\begin{center}
$W_{P_{\{i\}}}=\{w\in W_{G}\mid w(\omega _{i})=\omega _{i}\}$.
\end{center}

\noindent Via the embedding (3.3), the subgroup $W_{P_{\{i\}}}$ acts
naturally on $G/T$. The following lemma shows that every Weyl symmetry
arising from $W_{P_{\{i\}}}$ becomes trivial after passing to the base $%
G/P_{\{i\}}$.

\bigskip

\noindent \textbf{Lemma 3.1. }For each $w\in W_{P_{\{i\}}}$ the
diffeomorphism $w:G/T\rightarrow G/T$ descends to the identity on $%
G/P_{\{i\}}$, making the following diagram commutes:

\begin{center}
$%
\begin{array}{ccc}
G/T & \overset{w}{\rightarrow } & G/T \\ 
\pi \downarrow &  & \downarrow \pi \\ 
G/P_{\{i\}} & \overset{\mathrm{id}}{\rightarrow } & G/P_{\{i\}}%
\end{array}%
$.
\end{center}

\noindent \textbf{Proof. }Since $W_{G}=N_{G}(T)/T$, every element of $W_{G}$
admits a representative in $N_{G}(T)\subset G$. For a $w\in W_{G}$ choose a
representative $n_{w}\in N_{G}(T)$. Then the canonical right action of $w$
on $G/T$ is given by

\begin{center}
$w(gT)=gn_{w}T$, $g\in G$.
\end{center}

\noindent Likewise, since $W_{P_{\{i\}}}=N_{P_{\{i\}}}(T)/T$, every $w\in
W_{P_{\{i\}}}$ admits a representative

\begin{center}
$n_{w}\in N_{P_{\{i\}}}(T)\subset P_{\{i\}}$.
\end{center}

\noindent Consequently,

\begin{center}
$gn_{w}P_{\{i\}}=gP_{\{i\}},$ $g\in G$,
\end{center}

\noindent implying that the action of $w$ on $G/T$ descends via the map $\pi 
$ to the identity on $G/P_{\{i\}}$.\hfill $\square $

\subsection{The Dynkin symmetries on $G/T$}

It is a classical fact that every automorphism $\sigma \in \mathrm{Aut}%
(\Gamma _{G})$ of the Dynkin diagram $\Gamma _{G}$ extends uniquely to a Lie
algebra automorphism

\begin{center}
$\theta _{\sigma }:L(G)\rightarrow L(G)$
\end{center}

\noindent that preserves the Cartan subalgebra $L(T)$ (see, for example, 
\cite[\S 14]{Hu}). Since $G$ is simply connected, $\theta _{\sigma }$
integrates uniquely to a Lie group automorphism

\begin{center}
$\Theta _{\sigma }:G\rightarrow G$,
\end{center}

\noindent whose differential at the identity $e\in G$ is $\theta _{\sigma }$%
, and satisfying $\Theta _{\sigma }(T)=T$. Consequently, $\Theta _{\sigma }$
induces a diffeomorphism

\begin{center}
$\overline{\Theta }_{\sigma }:G/T\rightarrow G/T$.
\end{center}

\noindent Hence we obtain an injective homomorphism

\begin{enumerate}
\item[(3.4)] $\mathrm{Aut}(\Gamma _{G})\hookrightarrow \mathrm{Diff}(G/T)$, $%
\sigma \rightarrow \overline{\Theta }_{\sigma }$.
\end{enumerate}

\noindent The diffeomorphism $\overline{\Theta }_{\sigma }$ is called the 
\textsl{Dynkin symmetry} of $G/T$ associated with $\sigma $.

Lemma 3.1 shows that Weyl symmetries arising from the parabolic subgroup
become trivial after passing to the generalized Grassmannian. In contrast,
Dynkin symmetries fixing the fundamental weight $\omega _{{i}}$ descend to
genuine symmetries of $G/P_{\{i\}}$, which play a central role in the
statements and proofs of Theorems A and B.

\bigskip

\noindent \textbf{Lemma 3.2.} For every $\sigma \in \mathrm{Aut}(\Gamma
_{G},\omega _{{i}})$, the Lie group automorphism $\Theta _{\sigma }$
descends to a diffeomorphism

\begin{center}
$h_{\sigma }:G/P_{\{i\}}\rightarrow G/P_{\{i\}}$
\end{center}

\noindent making the following diagram commute:

\begin{center}
$%
\begin{array}{ccc}
G/T & \overset{\overline{\Theta }_{\sigma }}{\rightarrow } & G/T \\ 
\pi \downarrow &  & \downarrow \pi \\ 
G/P_{\{i\}} & \overset{h_{\sigma }}{\rightarrow } & G/P_{\{i\}}%
\end{array}%
$.
\end{center}

\noindent \textbf{Proof. }Since $\sigma \in \mathrm{Aut}(\Gamma _{G},\omega
_{{i}})$, the induced automorphism of $L(T)$ satisfies

\begin{center}
$\theta _{\sigma }(\omega _{{i}})=\sigma (\omega _{{i}})=\omega _{{i}}$.
\end{center}

\noindent For any $g\in P_{\{i\}}$, the definition of $P_{\{i\}}$ gives

\begin{center}
$g\cdot \exp (t\omega _{i})\cdot g^{-1}=\exp (t\omega _{i})$, $t\in \mathbb{R%
}$.
\end{center}

\noindent Applying $\Theta _{\sigma }$ to both sides and using

\begin{center}
$\Theta _{\sigma }(\exp U)=\exp \theta _{\sigma }(U)$,
\end{center}

\noindent we obtain

\begin{center}
$\Theta _{\sigma }(g)\exp (t\omega _{i})\Theta _{\sigma }(g)^{-1}=\exp
(t\omega _{i})$.
\end{center}

\noindent Thus $\Theta _{\sigma }(P_{\{i\}})\subseteq P_{\{i\}}$. Applying
the same argument to $\sigma ^{-1}$ gives the reverse inclusion, and hence

\begin{center}
$\Theta _{\sigma }(P_{\{i\}})=P_{\{i\}}$.
\end{center}

\noindent Therefore, $\Theta _{\sigma }$ descends to a diffeomorphism

\begin{center}
$h_{\sigma }:$ $G/P_{\{i\}}\rightarrow G/P_{\{i\}}$
\end{center}

\noindent and the commutative diagram follows immediately from the
definition of the quotient projection $\pi $.\hfill $\square $

\subsection{The stabilizer $\mathrm{Aut}(\Phi ,\protect\omega _{i})$}

Recall that the automorphism group $\mathrm{Aut}(\Phi )$ of the root system $%
\Phi $ admits the semidirect product decomposition

\begin{enumerate}
\item[(3.5)] $\mathrm{Aut}(\Phi )\cong W_{G}\rtimes \mathrm{Aut}(\Gamma
_{G}) $.
\end{enumerate}

\noindent See, for example, \cite[Ch.VI, \S 1.5, Prop.16]{Bou}. Combining
the Weyl symmetries (3.3) with the Dynkin symmetries (3.4) yields an
injective homomorphism

\begin{center}
$\mathrm{Aut}(\Phi )=W_{G}\rtimes \mathrm{Aut}(\Gamma _{G})\hookrightarrow 
\mathrm{Diff}(G/T)$.
\end{center}

\noindent Denote the resulting action by $\varphi \rightarrow h_{\varphi }$.
The following result is due to Papadima \cite[1.2 Theorem]{P}.

\bigskip

\noindent \textbf{Theorem 3.3. }The map $\varphi \rightarrow h_{\varphi
}^{\ast }$ is an anti-isomorphism

\begin{center}
$\mathrm{Aut}(\Phi )\rightarrow \mathrm{Aut}(H^{\ast }(G/T))$,
\end{center}

\noindent where the reversal of composition is due to the contravariance of
the pullback.\hfill $\square $

\bigskip

The canonical action of the automorphism group $\mathrm{Aut}(\Phi )$ on the
Cartan subalgebra $L(T)$ preserves the weight lattice $\Lambda \subset L(T)$%
. For this action of $\mathrm{Aut}(\Phi )$ on $\Lambda $, we prove the
following result, which provides the link needed to deduce Theorem B from
Theorem 3.3.

\bigskip

\noindent \textbf{Theorem 3.4.} For every fundamental weight $\omega _{i}\in
\Omega $, the stabilizer of $\omega _{i}$ in $\mathrm{Aut}(\Phi )$ admits
the semidirect product decomposition

\begin{enumerate}
\item[(3.6)] $\mathrm{Aut}(\Phi ,\omega _{i})\cong W_{P_{\left\{ i\right\}
}}\rtimes \mathrm{Aut}(\Gamma _{G},\omega _{i})$.
\end{enumerate}

\noindent \textbf{Proof.} Every element of $\mathrm{Aut}(\Phi )$ admits a
unique decomposition

\begin{center}
$\varphi =w\circ \sigma $, $w\in W_{G}$, $\sigma \in \mathrm{Aut}(\Gamma
_{G})$.
\end{center}

\noindent Suppose that $\varphi (\omega _{i})=\omega _{i}$. Then

\begin{center}
$\omega _{i}=w(\sigma (\omega _{i}))$.
\end{center}

\noindent Since $\sigma (\omega _{i})$ is again a fundamental weight, and
every Weyl orbit contains a unique dominant weight \cite[\S 10.1]{Hu}, it
follows that

\begin{center}
$\sigma (\omega _{i})=\omega _{i}$.
\end{center}

\noindent Consequently, $w(\omega _{i})=\omega _{i}$ and hence

\begin{center}
$w\in W_{P_{\left\{ i\right\} }}$, $\sigma \in \mathrm{Aut}(\Gamma
_{G},\omega _{i})$.
\end{center}

Moreover, since $W_{P_{\{i\}}}$ is normalized by $\mathrm{Aut}(\Gamma
_{G},\omega _{i})$

\begin{center}
$\sigma W_{P_{\{i\}}}\sigma ^{-1}=W_{P_{\{i\}}}$, $\sigma \in \mathrm{Aut}%
(\Gamma _{G},\omega _{i})$,
\end{center}

\noindent the uniqueness of the decomposition (3.5) implies

\begin{center}
$W_{P_{\left\{ i\right\} }}\cap \mathrm{Aut}(\Gamma _{G},\omega
_{i})=\left\{ \mathrm{id}\right\} $.
\end{center}

\noindent Therefore, $\mathrm{Aut}(\Phi ,\omega _{i})\cong W_{P_{\left\{
i\right\} }}\rtimes \mathrm{Aut}(\Gamma _{G},\omega _{i})$ as claimed.\hfill 
$\square $

\section{The basis theorem in Schubert calculus}

For a Lie group $G$ with a fixed maximal torus $T$, let $W_{G}\subset 
\mathrm{Aut}(L(T))$ denote its Weyl group, and let

\begin{center}
$l:W_{G}\rightarrow \mathbb{Z}_{\geq 0}$
\end{center}

\noindent be the\textsl{\ length function} on $W_{G}$. For every parabolic
subgroup $P\subset G$, the Weyl group $W_{P}$ is naturally identified with a
subgroup of $W_{G}$. Moreover, every right coset $wW_{P}$ contains a unique
element of minimal length. Accordingly, we identify the quotient $W_{G}/W_{P}
$ with the set of minimal-length representatives (see \cite[5.1]{BGG})

\begin{center}
$W(P,G):=\left\{ w\in W_{G}\mid l(w)\leq l(ww^{\prime }),w^{\prime }\in
W_{P}\right\} $,
\end{center}

\noindent Throughout, we regard $W(P,G)$ as a subset of $W_{G}$ via this
identification, and use the same notation $l$ for the restriction of the
length function to $W(P,G)$.

The following Schubert decomposition of $G/P$ was first established by
Ehresmann \cite{Eh} for the Grassmannians $G_{n,k}$ of $k$-dimensional
subspaces of $\mathbb{C}^{n}$, subsequently extended to the complete flag
manifolds $G/T$ by Chevalley \cite{Ch} and Demazure \cite{De}, and
established in full generality for arbitrary flag manifolds by
Bernstein--Gelfand--Gelfand \cite[\S 5]{BGG}.

\bigskip

\noindent \textbf{Lemma 4.1. }The flag manifold $G/P$\ admits a
decomposition into the Schubert subvarieties $X_{w}$ indexed by $W(P,G)$,

\begin{enumerate}
\item[(4.1)] $G/P=\underset{w\in W(P,G)}{\cup }X_{w}$, $\dim _{\mathbb{R}%
}X_{w}=2l(w)$,
\end{enumerate}

\noindent where $X_{w}$ is isomorphic to the closure of the Schubert cell
associated with $w\in W(P,G)$.\hfill $\square $

\bigskip 

Since the Schubert decomposition (4.1) consists entirely of the closures of
even-dimensional cells, the fundamental classes

\begin{center}
$[X_{w}]\in H_{2l(w)}(G/P)$, $\ w\in W(P,G)$,
\end{center}

\noindent form an additive basis of the integral homology $H_{\ast }(G/P)$.
Let $S_{w}\in H^{2l(w)}(G/P)$ denote the cohomology class Kronecker dual to $%
[X_{w}]$, characterized by

\begin{center}
$\left\langle S_{w},[X_{u}]\right\rangle =\delta _{w,u}$, $w,u\in W(P,G)$.
\end{center}

\noindent Lemma 4.1 immediately implies the following result, commonly known
as the \textsl{basis theorem} in Schubert calculus.

\bigskip

\noindent \textbf{Theorem 4.2.} The set of Schubert classes $\{S_{w},w\in
W(P,G)\}$ forms an additive basis of the integral cohomology $H^{\ast }(G/P)$%
.\hfill $\square $

\bigskip

The proofs of Theorems A and B rely on three consequences of the basis
theorem, which we establish in this section as Theorems 4.5, 4.7, and 4.9.
To provide the geometric foundation for these results, we begin by recalling
Bott--Samelson's desingularization of Schubert subvarieties.

\subsection{Bott-Samelson's desingularization of Schubert varieties}

Retaining the notation of Section 3, let $G$ be a simply connected simple
Lie group with a fixed maximal torus $T$, and let $x_{0}\in L(T)$ be a
regular point. Denote by $\mathcal{F}(x_{0})$ the closure of the Weyl
chamber containing $x_{0}$, and let

\begin{center}
$L(x_{0})=\{L_{1},\cdots ,L_{n}\}\subset L(T)$
\end{center}

\noindent be the collection of walls of $\mathcal{F}(x_{0})$ (cf. (3.4)).
The following facts are standard. For each $1\leq i\leq n$ let $s_{i}\in 
\mathrm{Aut}(L(T))$ be the reflection across the hyperplane $L_{i}$, and let 
$K_{{i}}\subset G$ denote the centralizer of the subtorus $\exp (L_{{i}})$.

\bigskip

\noindent \textbf{Lemma 4.3. }(i) The Weyl group $W_{G}$ is naturally
isomorphic to the subgroup of $\mathrm{Aut}(L(T))$ generated by the
reflections $s_{i}$, $1\leq i\leq n$.

(ii)\ Moreover, $T$ is a maximal torus of $K_{{i}}$, and the homogeneous
space $K_{{i}}/T$ is diffeomorphic to the $2$-sphere.\hfill $\square $

\bigskip

Let $P\subset G$ be a parabolic subgroup. By Lemma 4.3(i), every left coset
of $W_{G}/W_{P}$ contains a unique representative $w$ of minimal length. Fix
a reduced expression

\begin{enumerate}
\item[(4.2)] $w=s_{i_{1}}\circ \cdots \circ s_{i_{k}}\in W(P,G)$,
\end{enumerate}

\noindent where $k=l(w)$ and $1\leq i_{j}\leq n$. Associated with this
reduced decomposition, define

\begin{center}
$\varphi _{w}:K_{w}:=K_{i_{1}}\times \cdots \times K_{i_{k}}\rightarrow G$
\end{center}

\noindent by

\begin{center}
$\ \varphi _{w}(g_{1},\cdots ,g_{k})=g_{1}\cdots g_{k}$.
\end{center}

\noindent Since $T\subset K_{{i}}$ for every $i$ by Lemma 4.3(ii), the
product $T\times \cdots \times T$ ($k$-copies) acts freely on $K_{w}$ by

\begin{center}
$(g_{{1}},\ldots ,g_{{k}})(t_{{1}},\ldots ,t_{{k}})=(g_{{1}}t_{{1}},t_{{1}%
}^{-1}g_{{2}}t_{{2}},\ldots ,t_{{k-1}}^{-1}g_{{k}}t_{{k}})$.
\end{center}

\noindent Let $\Gamma _{w}:=$ $K_{w}/T^{k}$ denote the quotient manifold
endowed with the orientation induced from the canonical orientations of $%
K_{w}$ and $T^{k}$. We write

\begin{center}
$[g_{{1}},\ldots ,g_{{k}}]\in \Gamma _{\overline{w}}$
\end{center}

\noindent for the point represented by $(g_{{1}},\ldots ,g_{{k}})$. The
associated \textsl{Bott-Samelson cycle} is the map

\begin{center}
$\overline{\varphi }_{w,P}:$ $\Gamma _{w}\rightarrow G/P$
\end{center}

\noindent defined by

\begin{center}
$\overline{\varphi }_{w,P}([g_{{1}},\ldots ,g_{{k}}])=g_{{1}}\cdots g_{{k}}P$%
.
\end{center}

By Lemma 1.2, we may assume that $P=P_{I}$ for some subset $I\subseteq
\{1,\cdots ,n\}$. Let

\begin{center}
$\pi _{I}:G/T\rightarrow G/P_{I}$
\end{center}

\noindent be the canonical projection. Observe that when $I=\{1,\cdots ,n\}$,

\begin{center}
$\pi _{I}=\mathrm{id}:G/T\rightarrow G/T$.
\end{center}

\noindent The following result was proved by Hansen \cite{Ha} for the case $%
P=T$. The general case follows immediately by composing the Bott-Samelson
cycle with the canonical projection $\pi _{I}$ (see \cite[Lemma 6.1.]{D1}).

\bigskip

\noindent \textbf{Lemma 4.4.} For every $w\in W(P,G)$, the Bott--Samelson
cycle

\begin{center}
$\overline{\varphi }_{w,P}:$ $\Gamma _{w}\rightarrow G/P$
\end{center}

\noindent is a proper surjective map of degree one onto the Schubert variety 
$X_{w}$, and

\begin{enumerate}
\item[(4.3)] $\overline{\varphi }_{w,P\ast }[\Gamma _{w}]=[X_{w}]\in
H_{2l(w)}(G/P)$.\hfill $\square $
\end{enumerate}

Because of the relation (4.3), the Bott-Samelson cycle $\overline{\varphi }%
_{w,P}$ is also called a \textsl{Bott--Samelson desingularization} (or 
\textsl{resolution}) of the Schubert subvariety $X_{w}$.

\subsection{The action of the Dynkin symmetries on $H^{\ast }(G/T)$}

Continuing the discussion from Section 3.3, let $\overline{\Theta }_{\sigma
}:G/T\rightarrow G/T$ be a Dynkin symmetry associated with an automorphism $%
\sigma \in \mathrm{Aut}(\Gamma _{G})$. Since $\sigma $ permutes the set $%
\Delta =\{\alpha _{1},\cdots ,\alpha _{n}\}$ of simple roots, the induced
action $\theta _{\sigma }$ on the Cartan subalgebra $L(T)$ permutes the walls

\begin{center}
$L(x_{0})=\{L_{1},\cdots ,L_{n}\}$
\end{center}

\noindent of the Weyl chamber $\mathcal{F}(x_{0})$ according to the rule

\begin{center}
$L_{i}\rightarrow L_{\sigma (i)}$ if $\sigma (\alpha _{i})=\alpha _{\sigma
(i)}$, for $1\leq i\leq n$.
\end{center}

\noindent The action of the Dynkin symmetry $\overline{\Theta }_{\sigma }$
on the Schubert basis of $H^{\ast }(G/T)$ is described by the following
result.

\bigskip

\noindent \textbf{Theorem 4.5. }Let $\sigma \in \mathrm{Aut}(\Gamma _{G})$
be an automorphism, and let $w=s_{i_{1}}\cdots s_{i_{k}}$ be a reduced
decomposition of $w\in W_{G}$. Define

\begin{center}
$\sigma (w)=s_{\sigma (i_{1})}\cdots s_{\sigma (i_{k})}$.
\end{center}

\noindent Then

\begin{enumerate}
\item[(4.4)] $\overline{\Theta }_{\sigma }^{\ast }(S_{w})=S_{\sigma (w)}$.
\end{enumerate}

\noindent \textbf{Proof. }Since $\sigma $ preserves the Coxeter relations
among the simple reflections $s_{i}$'s, it induces an automorphism of the
Weyl group $W_{G}$. Hence the image of every reduced expression is again
reduced. Therefore,

\begin{center}
$l(\sigma (w))=l(w)$.
\end{center}

Since $\theta _{\sigma }(L_{i})=L_{\sigma (i)}$, the automorphism $\Theta
_{\sigma }$ sends the subtorus $\exp (L_{{i}})$ onto $\exp (L_{\sigma ({i)}%
}) $. Hence

\begin{enumerate}
\item[(4.5)] $\Theta _{\sigma }(K_{i})=K_{\sigma (i)}$.
\end{enumerate}

\noindent It follows that the maps $\varphi _{w}$ and $\Theta _{\sigma }$
satisfy the commutative diagram

\begin{center}
$%
\begin{array}{ccc}
K_{w}=K_{i_{1}}\times \cdots \times K_{i_{k}} & \overset{\varphi _{w}}{%
\rightarrow } & G \\ 
\Theta \downarrow &  & \downarrow \Theta _{\sigma } \\ 
K_{\sigma (w)}=K_{\sigma (i_{1})}\times \cdots \times K_{\sigma (i_{k})} & 
\overset{\varphi _{\sigma (w)}}{\rightarrow } & G%
\end{array}%
$,
\end{center}

\noindent where the left vertical map $\Theta $ is the product of the
diffeomorphisms $K_{i}\cong K_{\sigma (i)}$ arising from (4.5). By Lemma
4.4, passing to the quotient manifolds gives the commutative diagram,
showing that

\begin{center}
$\overline{\Theta }_{\sigma }(X_{w})=X_{\sigma (w)}$.
\end{center}

\noindent Since the Schubert classes are defined as the Kronecker duals of
the fundamental classes of the Schubert varieties, this identity immediately
yields (4.4) by naturality of the Kronecker pairing.$\hfill \square $

\subsection{Borel--Hirzebruch isomorphism $\Lambda \protect\cong H^{2}(G/T)$}

Let $\Lambda $ be the weight lattice of $L(T)$, and let $\varpi _{i}\in $ $%
H^{2}(G/T)$ denote the Schubert class corresponding to the simple reflection 
$s_{i}\in W_{G}$. Since $\Lambda $ is freely generated by the set $\Omega
=\{\omega _{1},\cdots ,\omega _{n}\}$ of fundamental dominant weights, we
can introduce the additive isomorphism $h$

\begin{center}
$h:\Lambda \rightarrow H^{2}(G/T)$
\end{center}

\noindent by $h(\omega _{i})=$ $\varpi _{i}$, $1\leq i\leq n$. Recall that
the automorphism group $\mathrm{Aut}(\Phi )$ acts canonically on $\Lambda $,
while its action on $H^{2}(G/T)$ is induced by the Dynkin and Weyl
symmetries described in Section 3.

\bigskip

\noindent \textbf{Lemma 4.6. }The map $h$ is an isomorphism of $\mathrm{Aut}%
(\Phi )$-modules.

\bigskip

\noindent \textbf{Proof. }Let\textbf{\ }$A=(a_{i,j})\in M_{n}(\mathbb{Z})$
be the Cartan matrix of $G$. The action of the simple reflection $s_{i}\in
W_{G}$ on $\Lambda $ on the fundamental weights is given by

\begin{center}
$s_{i}(\omega _{k})=\left\{ 
\begin{tabular}{l}
$\omega _{k}$ if $i\neq k$; \\ 
$\omega _{i}-(a_{i,1}\omega _{1}+\cdots +a_{i,n}\omega _{n})$ if $i=k$.%
\end{tabular}%
\right. $
\end{center}

\noindent On the other hand, the BGG formula \cite[3.14 Theorem(iii)]{BGG}
gives the same action of $s_{i}$ on the Schubert basis $\{\varpi _{k}\}$.
Hence $h$ is $W_{G}$-equivariant.

Furthermore, by formula (4.4), $h$ is also an isomorphism of $\mathrm{Aut}%
(\Gamma _{G})$-modules. Since both module structures extend to the
semidirect product $W_{G}\rtimes \mathrm{Aut}(\Gamma _{G})$, the map $h$ is $%
\mathrm{Aut}(\Phi )$-equivariant.$\hfill \square $

\bigskip

Combining Theorem 3.3, Theorem 3.4, and Lemma 4.6, we obtain the following
description of the stabilizer of $\varpi _{i}\in H^{2}(G/T)$ in $\mathrm{Aut}%
(H^{\ast }(G/T))$, which is the key algebraic ingredient in the proof of
Theorem B.

\bigskip

\noindent \textbf{Theorem 4.7. }For a Schubert basis element $\varpi _{i}\in
H^{2}(G/T)$, there is an anti-isomorphism

\begin{enumerate}
\item[(4.5)] $\mathrm{Aut}(H^{\ast }(G/T),\varpi _{i})\cong W_{P_{\left\{
i\right\} }}\rtimes \mathrm{Aut}(\Gamma _{G},\omega _{i})$.\hfill $\square $
\end{enumerate}

\bigskip

\noindent \textbf{Remark 4.8. }The map $h:\Lambda \rightarrow H^{2}(G/T)$
was introduced by Borel-Hirzebruch \cite[\S 10.2]{BH} as follows. Every
weight $\omega \in \Lambda $ determines a character $T\rightarrow S^{1}$,
and hence a principal $S^{1}$-bundle over $G/T$. The first Chern class of
this bundle is precisely

\begin{center}
$h(\omega )\in H^{2}(G/T)$.
\end{center}

\noindent It is well known that $h$ is an isomorphism of $W_{G}$-modules.

Lemma 4.6 shows that the Borel--Hirzebruch isomorphism enjoys two additional
properties that deserve emphasis:

(i) it is equivariant with respect to the full automorphism group $\mathrm{%
Aut}(\Phi )$;

(ii) it identifies the fundamental dominant weights $\Omega =\{\omega
_{1},\cdots ,\omega _{n}\}$ with the Schubert basis of $H^{2}(G/T)$.

\noindent These two properties together identify the stabilizer of a
Schubert class $\varpi _{i}$ with the stabilizer of the corresponding
fundamental weight $\omega _{i}$, yielding the isomorphism (4.5).\hfill $%
\square $

\subsection{The Leray-Hirsch property of $\protect\pi :G/T\rightarrow G/P$}

For a subset $I\subseteq \{1,\cdots ,n\}$ let $P=P_{I}$ be the corresponding
parabolic subgroup of $G$, and consider the fiber bundle induced by the
inclusions $T\subset P$ $\subset G$

\begin{center}
$P/T\overset{i}{\hookrightarrow }G/T\overset{\pi }{\rightarrow }G/P$.
\end{center}

\noindent By Bott--Samelson's description of Schubert varieties (Lemma 4.4),
the induced maps $\pi ^{\ast }$ and $i^{\ast }$ are compatible with the
Schubert bases of the three flag manifolds

\begin{center}
$P/T$, $G/P$ and $G/T$.
\end{center}

\noindent More precisely, one has:

\begin{quote}
(i) Under the inclusion $W_{P}\subset W_{G}$,\ the classes $S_{w}$, $w\in
W_{P}$, restrict to the Schubert basis of $H^{\ast }(P/T)$, and hence form a
Leray--Hirsch basis for $\pi $.

(ii) Under the inclusion $W(P,G)\subset W_{G}$,\ the map $\pi ^{\ast }$
identifies the Schubert basis $\{S_{w}\}_{w\in W(P;G)}$\ of $H^{\ast }(G/P)$%
\ with the corresponding subset of the Schubert basis of $H^{\ast }(G/T)$.
\end{quote}

\noindent Consequently, the bundle

\begin{center}
$\pi :G/T\rightarrow G/P$
\end{center}

\noindent satisfies the Leray--Hirsch property over the integers. It follows
from the Leray--Hirsch theorem that the cohomology $H^{\ast }(G/T)$ is a
free $H^{\ast }(G/P)$-module with the basis $\{S_{w}\}_{w\in W_{P}}$, namely

\begin{enumerate}
\item[(4.6)] $H^{\ast }(G/T)=H^{\ast }(G/P)\{S_{w}\}_{w\in W_{P}}$.
\end{enumerate}

Now let $f:$ $G/P\rightarrow $ $G/P$ be a self-map, and consider its induced
bundle $\pi _{f}$

\begin{enumerate}
\item[(4.7)] $%
\begin{array}{ccc}
E_{f} & \overset{\widetilde{f}}{\rightarrow } & G/T \\ 
\pi _{f}\downarrow  &  & \downarrow \pi  \\ 
G/P & \overset{f}{\rightarrow } & G/P%
\end{array}%
$,
\end{enumerate}

\noindent where $\widetilde{f}$ denotes the canonical bundle map covering $f$%
. Since $\widetilde{f}$ restricts to identity on each fiber, we may use the
same notation $S_{w}\in H^{\ast }(E_{f})$ for the pullback $\widetilde{%
f^{\ast }}(S_{w})$.

\bigskip

\noindent \textbf{Theorem 4.9. }The integral cohomology ring of $E_{f}$ is a
free $H^{\ast }(G/P)$-module with the basis $\{S_{w}\}_{w\in W_{P}}$

\begin{center}
$H^{\ast }(E_{f})=H^{\ast }(G/P)\{S_{w}\}_{w\in W_{P}}$.
\end{center}

\noindent With respect to this decomposition, the induced map $\widetilde{f}%
^{\ast }$ on cohomology is given by

\begin{center}
$\widetilde{f}^{\ast }(x\cup S_{w})=f^{\ast }(x)\cup S_{w}$, $x\in H^{\ast
}(G/P)$, $w\in W_{P}$.
\end{center}

In particular, if $f$ is a rational homotopy equivalence, then the bundle
map $\widetilde{f}$ is also a rational homotopy equivalence.

\bigskip

\noindent \textbf{Proof.} Since pullbacks preserve fibers, the fiber of $\pi
_{f}$ is again $P/T$. Moreover, the pullback of a Leray--Hirsch bundle is
again a Leray--Hirsch bundle. Applying (4.6) to the pullback bundle
immediately gives

\begin{center}
$H^{\ast }(E_{f})=H^{\ast }(G/P)\{s_{w}\}_{w\in W_{P}}$.
\end{center}

The formula for $\widetilde{f}^{\ast }$ follows from the commutative diagram

\begin{center}
$\widetilde{f}^{\ast }\circ \pi ^{\ast }=\pi _{f}^{\ast }\circ f^{\ast }$,
\end{center}

\noindent together with the fact that the canonical bundle map $\widetilde{f}
$ identifies each fiber of $\pi _{f}$ with the corresponding fiber of $\pi $
by the identity map.

If $f$ is a rational homotopy equivalence, then

\begin{center}
$f^{\ast }:$ $H^{\ast }(G/P;\mathbb{Q})\rightarrow H^{\ast }(G/P;\mathbb{Q})$
\end{center}

\noindent is an isomorphism. The preceding decomposition shows that $%
\widetilde{f}^{\ast }$ is likewise an isomorphism on rational cohomology.
Hence the bundle map

\begin{center}
$\widetilde{f}:$ $E_{f}\rightarrow G/T$
\end{center}

\noindent is a rational homotopy equivalence.\hfill $\square $

\bigskip 

For a smooth manifold $M$,\ let $M_{0}$ denote its rationalization.
Likewise, for a smooth map $p:E\rightarrow B$, let

\begin{center}
$p_{0}:E_{0}\rightarrow B_{0}$
\end{center}

\noindent denote the rationalization of $p$. When $p$ is a smooth fiber
bundle with simply connected fiber, its rationalization $p_{0}$ need not be
a fiber bundle, but is a Serre fibration with a connected nilpotent fiber of
finite rational type. 

Applying rationalization to the pullback diagram (4.7), we obtain

\begin{center}
$%
\begin{array}{ccc}
(E_{f})_{0} & \overset{\widetilde{f}_{0}}{\rightarrow } & (G/T)_{0} \\ 
(\pi _{f})_{0}\downarrow  &  & \downarrow \pi _{0} \\ 
(G/P)_{0} & \overset{f_{0}}{\rightarrow } & (G/P)_{0}%
\end{array}%
$,
\end{center}

\noindent where both rationalized maps $(\pi _{f})_{0}$ and $\pi _{0}$ are
Serre fibrations with connected nilpotent fibers of finite rational type.

We shall use the standard fact that the pullback of such a fibration along a
homotopy equivalence of the base is fiber-homotopy equivalent to the
original fibration after rationalization; see \cite{HT}. Thus we obtain the
following.

\bigskip 

\noindent \textbf{Corollary 4.10. }If $f:G/P\rightarrow G/P$ is a rational
self-homotopy equivalence, there exists a rational homotopy equivalence

\begin{center}
$\varepsilon :(G/T)_{0}\rightarrow (E_{f})_{0}$
\end{center}

\noindent over the identity of $(G/P)_{0}$, making the following diagram
commute

\begin{center}
$%
\begin{array}{ccc}
(G/T)_{0} & \overset{\varepsilon }{\rightarrow } & (E_{f})_{0} \\ 
\pi _{0}\downarrow  &  & \downarrow (\pi _{f})_{0}\quad  \\ 
(G/P)_{0} & \overset{\mathrm{id}}{\rightarrow } & (G/P)_{0}%
\end{array}%
$.
\end{center}

\noindent \textbf{Proof. }Since $P/T$ is simply connected, both $(\pi
_{f})_{0}$ and $\pi _{0}$ are Serre fibrations with connected nilpotent
fibers of finite rational type. Moreover, $(\pi _{f})_{0}$ is the pullback
of $\pi _{0}$ along the homotopy equivalence

\begin{center}
$f_{0}:(G/P)_{0}\rightarrow (G/P)_{0}$. 
\end{center}

\noindent By the rational fiber-homotopy invariance of such pullback
fibrations \cite{HT}, the fibrations $(\pi _{f})_{0}$ and $\pi _{0}$ are
fiber-homotopy equivalent over $(G/P)_{0}$. Hence there exists a rational
homotopy equivalence

\begin{center}
$\varepsilon :(G/T)_{0}\rightarrow (E_{f})_{0}$
\end{center}

\noindent over the identity of $(G/P)_{0}$, as asserted.\hfill $\square $

\bigskip 

\noindent \textbf{Remark 4.11. }In his fifteenth problem, Hilbert called for
a rigorous foundation of Schubert's enumerative calculus. In modern terms,
Weil \cite[p.331]{W} reformulated this problem as the determination of the
integral cohomology ring of a flag manifold. In view of the basis theorem
(Theorem 4.2), this amounts to expressing products of Schubert classes as
integral linear combinations of the Schubert basis. Schubert himself
referred to this as the characteristic problem \cite{Sch1,Sch2}, and it has
long been regarded as a central problem of Schubert calculus (see, for
example, \cite{Eh,Se,Wa}). Solutions to both the characteristic problem and
Weil's formulation, together with historical accounts of their development,
can be found in the survey articles \cite{DZ3,DZ4}.

The present paper provides another application of the basis theorem and its
geometric consequences: the classification of self-maps of generalized
Grassmannians.\hfill $\square $

\section{The proofs of Theorems A and B}

Let $P\subset G$ be a parabolic subgroup. Every element $w\in W(P,G)$ admits
a reduced decomposition

\begin{center}
$w=s_{i_{1}}\circ \cdots \circ s_{i_{k}}$,
\end{center}

\noindent where $k=l(w)$ and $1\leq i_{j}\leq n$. Accordingly, we write

\begin{center}
$w=s_{J}$, $J=(i_{{1}},\ldots ,i_{{k}})$.
\end{center}

\noindent Such reduced decomposition of a given $w$ is generally not unique.
To remove this ambiguity, let $D(w)$ denote the set of all reduced
decompositions of $w$. Since $D(w)$ is finite, it has a unique
lexicographically minimal element. We call this \textsl{the
lexicographically minimal reduced decomposition}, or simply \textsl{the
minimal decomposition}, of $w$.

Furthermore, for each integer $r\geq 1$, set

\begin{center}
$W^{r}(P,G):=\{w\in W(P,G)\mid l(w)=r\}$, $\beta _{r}:=\left\vert
W^{r}(P,G)\right\vert $.
\end{center}

\noindent The lexicographic ordering of the minimal decompositions induces a
total order on $W^{r}(P,G)$. We therefore enumerate its elements as

\begin{enumerate}
\item[(5.1)] $W^{r}(P,G)=\{w_{{r,1}},\cdots ,w_{{r,}\beta _{r}}\}$,
\end{enumerate}

\noindent where $w_{{r,i}}$ denotes the $i$-th element in this ordering.

In our work on computing the Chow rings of generalized Grassmannians \cite%
{DZ1}, we implemented the program \textsl{Decomposition}, whose output is
summarized as follows.

\bigskip

\noindent \textbf{Algorithm 5.1 (Decomposition)}.

\begin{quote}
\textbf{Input:} \textsl{The Cartan matrix }$A=(a_{{ij}})_{{n\times n}}$ 
\textsl{of }$G$\textsl{, together with a subset }$I\subset \{1,\ldots ,n\}$%
\textsl{.}

\textbf{Output: }\textsl{The ordered sets }$W^{r}(P_{I},G)$\textsl{\ (}$%
r\geq 1)$\textsl{, with each element represented by its minimal
decomposition and indexed as in (5.1)}.\hfill $\square $
\end{quote}

\noindent

\noindent Examples of the output of Decomposition can be found in \cite[%
1.1--7.1]{DZ5}.

The proof of Theorem A relies on Algorithm 5.1, which provides an effective
ordering of $W(P_{\{i\}},G)$, first by length and then lexicographically.
This ordering reduces the detection of a nontrivial Dynkin symmetry to a
finite and effective search: one inspects the Schubert classes in increasing
order until finding the first $w\in W(P_{\{i\}},G)$ such that $\sigma
(w)\neq w$, where $\sigma \in \mathrm{Aut}(\Gamma _{G},\omega _{{i}})$ is
nontrivial. In particular, the algorithm provides a systematic method for
identifying a Schubert class of minimal dimension on which a given
nontrivial Dynkin symmetry acts nontrivially.

\bigskip

\noindent \textbf{Proof of Theorem A. }Let $(G,i)$ be as in Lemma 1.4. For a 
$0\neq \sigma \in \mathrm{Aut}(\Gamma _{G},\omega _{{i}})$ let

\begin{center}
$h_{\sigma }:G/P_{\{i\}}\rightarrow G/P_{\{i\}}$
\end{center}

\noindent be the associated Dynkin symmetry. By Theorem 4.5, it suffices to
find a Schubert class $S_{w}$ such that $\sigma (w)\neq w$. Algorithm 5.1
provides an effective way to locate such a class by examining $W(P_{\{i\}},G)
$ in increasing order of length and lexicographically within each length.
Lemma 1.4 determines all nontrivial stabilizers $\mathrm{Aut}(\Gamma
_{G},\omega _{{i}})$, so we distinguish two cases accordingly.

\textbf{Case 1.} If $(G,i)\neq (D_{4},2)$, then $\mathrm{Aut}(\Gamma
_{G},\omega _{{i}})\cong \mathbb{Z}_{2}$ by Lemma 1.4. Let $\sigma $ denote
the generator of this group. 

In each case of $(G,i)$, the second column of the following table records
the first element in the prescribed ordering for which $\sigma (w)\neq w$,
as shown in the third column. The corresponding Schubert class $S_{w}$
therefore provides an explicit detector in the lowest possible dimension.

\begin{center}
\begin{tabular}{l||l|l|}
\hline
$(G,i)$ & $w\in W(P_{\{i\}},G)$ & $\sigma (w)$ \\ \hline\hline
$(A_{2k+1},k+1)$ & $w_{2,1}=s_{\{k,k+1\}}$ & $w_{2,2}=s_{\{k+2,k+1\}}$ \\ 
\hline
$(E_{6},2)$ & $w_{3,1}=s_{\{3,4,2\}}$ & $w_{3,2}=s_{\{5,4,2\}}$ \\ \hline
$(E_{6},4)$ & $w_{2,2}=s_{\{3,4\}}$ & $w_{2,3}=s_{\{5,4\}}$ \\ \hline
$(D_{n},i)$, $n>4,1\leq i\leq n-2$ & $w_{n-i,1}=s_{\{n-1,n-2,\cdots
,i+1,i\}} $ & $w_{n-i,2}=s_{\{n,n-2,\cdots ,i+1,i\}}$ \\ \hline
$(D_{4},1)$ & $w_{3,1}=s_{\{3,2,1\}}$ & $w_{3,2}=s_{\{4,2,1\}}$ \\ \hline
$(D_{4},3)$ & $w_{3,1}=s_{\{1,2,3\}}$ & $w_{3,2}=s_{\{4,2,3\}}$ \\ \hline
$(D_{4},4)$ & $w_{3,1}=s_{\{1,2,4\}}$ & $w_{3,2}=s_{\{3,2,4\}}$ \\ \hline
\end{tabular}%
.
\end{center}

\textbf{Case 2.} If $(G,i)=(D_{4},2)$. Then $\mathrm{Aut}(\Gamma
_{D_{4}},\omega _{{i}})\cong S_{3}$ Lemma 1.4.

By \textsl{Decomposition} the set $W^{2}(P_{\{2\}},D_{4})$ consists of the
three elements

\begin{center}
$w_{2,1}=s_{\{1,2\}}$, $w_{2,2}=s_{\{3,2\}}$, $w_{2,3}=s_{\{4,2\}}$.
\end{center}

\noindent Moreover, we may choose generators $\sigma $, $\tau \in S_{3}$ of
orders $2$ and $3$, respectively, such that their actions on the set $%
W^{2}(P_{\{2\}},D_{4})$ are given by

\begin{center}
$\left\{ w_{2,1},w_{2,2},w_{2,3}\right\} \overset{\sigma }{\rightarrow }%
\left\{ w_{2,1},w_{2,3},w_{2,2}\right\} $,

$\left\{ w_{2,1},w_{2,2},w_{2,3}\right\} \overset{\tau }{\rightarrow }%
\left\{ w_{2,3},w_{2,1},w_{2,2}\right\} $,
\end{center}

\noindent respectively. In particular, both $\sigma $ and $\tau $ act
nontrivially on the Schubert basis, and hence

\begin{center}
$h_{\sigma }^{\ast }$, $h_{\tau }^{\ast }\neq \mathrm{Id}$.
\end{center}

\noindent This completes the proof of Theorem A.\hfill $\square $

\bigskip

Turning to the proof of Theorem B, let $G$ be a simply connected simple Lie
group with a fixed maximal torus $T$, and let $P_{\{i)}$ be the parabolic
subgroup associated with a fundamental dominant weight $\omega _{i}\in
\Omega $. Consider the fiber bundle

\begin{center}
$P_{\{i)}/T\overset{i}{\hookrightarrow }G/T\overset{\pi }{\rightarrow }%
G/P_{\{i)}$
\end{center}

\noindent induced by the inclusions $T\subset $ $P_{\{i)}\subset G$. In the
notation of Section 4.3, let

\begin{center}
$\{\varpi _{1},\cdots ,\varpi _{n}\}\subset H^{2}(G/T)$
\end{center}

\noindent denote the Schubert basis. As in Section 1, set $X=G/P_{\{i)}$.
Since the K\"{a}hler class

\begin{center}
$\kappa _{X}\in H^{2}(X)\cong \mathbb{Z}$
\end{center}

\noindent is precisely the Schubert basis element $S_{s_{i}}\in H^{2}(X)$,
the Leray--Hirsch property of $\pi $ stated in Section 4.4(ii) gives the
following.

\bigskip

\noindent \textbf{Lemma 5.2.} The induced homomorphism $\pi ^{\ast
}:H^{2}(X)\rightarrow H^{2}(G/T)$ satisfies

\begin{enumerate}
\item[(5.2)] $\pi ^{\ast }(\kappa _{X})=\varpi _{i}$.\hfill $\square $
\end{enumerate}

Since the ring $H^{\ast }(G/T)$ is torsion free by Theorem 4.2, the
inclusion $\mathbb{Z}\hookrightarrow \mathbb{Q}$ induces a group monomorphism

\begin{center}
$j:\mathrm{Aut}(H^{\ast }(G/T))\hookrightarrow \mathrm{Aut}(H^{\ast }(G/T;%
\mathbb{Q}))$.
\end{center}

\noindent The following result is essentially contained in the proofs of
Theorems 1.1 and 1.2 of Papadima \cite{P}. We state it here in the form
needed below.

\bigskip 

\noindent \textbf{Lemma 5.3.} For any nonzero element $\omega \in H^{2}(G/T)$%
, the monomorphism $j$ restricts to an isomorphism of stabilizers

\begin{enumerate}
\item[(5.3) ] $\mathrm{Aut}(H^{\ast }(G/T),\omega )\rightarrow \mathrm{Aut}%
(H^{\ast }(G/T;\mathbb{Q}),\omega )$.\hfill $\square $
\end{enumerate}

\noindent \textbf{Proof of Theorem B. }Let $X=G/P_{\{i)}$, and let $%
f:X\rightarrow X$ be a self-map with $\deg (f)\neq 0$. By Corollary 2.2, $f$
is a rational homotopy equivalence. Consider the pullback bundle

\begin{center}
$%
\begin{array}{ccc}
E_{f} & \overset{\widetilde{f}}{\rightarrow } & G/T \\ 
\pi _{f}\downarrow &  & \downarrow \pi \\ 
G/P & \overset{f}{\rightarrow } & G/P%
\end{array}%
$,
\end{center}

\noindent where the horizontal map $\widetilde{f}$ is the canonical bundle
map covering $f$. By Theorem 4.9, the rationalization

\begin{center}
$\widetilde{f}_{0}:(E_{f})_{0}\rightarrow (G/T)_{0}$ 
\end{center}

\noindent is a homotopy equivalence. By Corollary 4.10, there exists a
rational homotopy equivalence

\begin{center}
$\varepsilon :(G/T)_{0}\rightarrow (E_{f})_{0}$
\end{center}

\noindent over the identity of $X_{0}$. Define

\begin{center}
$F:=\varepsilon \circ \widetilde{f}_{0}:(G/T)_{0}\rightarrow (G/T)_{0}$.
\end{center}

\noindent Then $F$ is a rational self-homotopy equivalence covering $f_{0}$:

\begin{enumerate}
\item[(5.4)] $%
\begin{array}{ccc}
(G/T)_{0} & \overset{F}{\rightarrow } & (G/T)_{0} \\ 
\pi _{0}\downarrow  &  & \downarrow \pi _{0} \\ 
X_{0} & \overset{f_{0}}{\rightarrow } & X_{0}%
\end{array}%
$, i.e. $\pi _{0}\circ F=f_{0}\circ \pi _{0}$.
\end{enumerate}

\noindent In particular, 

\begin{center}
$F^{\ast }\in \mathrm{Aut}(H^{\ast }(G/T;\mathbb{Q}))$.
\end{center}

Via the canonical isomorphism

\begin{center}
$H^{2}((G/T)_{0};\mathbb{Q})\cong H^{2}(G/T;\mathbb{Q})$,
\end{center}

\noindent we continue to denote by $\{\varpi _{1},\cdots ,\varpi _{n}\}$ the
Schubert basis of $H^{2}((G/T)_{0};\mathbb{Q})$. Under this identification,
(5.2) gives

\begin{center}
$\pi _{0}^{\ast }(\kappa _{X})=\varpi _{i}$
\end{center}

\noindent Applying the cohomology functor to (5.4), we obtain

\begin{center}
$F^{\ast }(\varpi _{i})=F^{\ast }(\pi _{0}^{\ast }(\kappa _{X}))=\pi
_{0}^{\ast }(f_{0}^{\ast }(\kappa _{X}))=\deg (f)\cdot \varpi _{i}$.
\end{center}

\noindent In particular,

\begin{enumerate}
\item[(5.5)] $F^{\ast }(\varpi _{i})=\lambda \cdot \varpi _{i}$, $\lambda :=$
$\deg (f)\in \mathbb{Q}^{\times }$.
\end{enumerate}

Consider the automorphism

\begin{center}
$A:=\psi ^{\lambda ^{-1}}\circ F^{\ast }\in \mathrm{Aut}(H^{\ast }(G/T;%
\mathbb{Q}))$,
\end{center}

\noindent where, for $a\in \mathbb{Q}^{\times }$, the rational Adams
operator $\psi ^{a}$ is defined by

\begin{center}
$\psi ^{a}(x)=a^{r}x$, $x\in H^{2r}(G/T;\mathbb{Q})$.
\end{center}

\noindent By (5.5), 

\begin{center}
$A(\varpi _{i})=\psi ^{\lambda ^{-1}}(\lambda \cdot \varpi _{i})=\varpi _{i}$
\end{center}

\noindent Hence

\begin{center}
$A\in \mathrm{Aut}(H^{\ast }(G/T;\mathbb{Q}),\varpi _{i})$.
\end{center}

\noindent It follows from Lemma 5.3 that

\begin{center}
$A\in \mathrm{Aut}(H^{\ast }(G/T),\varpi _{i})$.
\end{center}

\noindent By Theorem 4.7, there exist $w\in W_{P_{\left\{ i\right\} }}$ and $%
\sigma \in \mathrm{Aut}(\Gamma _{G},\omega _{i})$ such that

\begin{center}
$A=\overline{\Theta }_{\sigma }^{\ast }\circ w^{\ast }$
\end{center}

By Lemma 3.1, $w$ is a fiberwise self-diffeomorphism of $G/T$ over the
identity of $X$. Replacing $\varepsilon $ by $w^{-1}\circ \varepsilon $
therefore preserves the condition that $\varepsilon $ is over $X_{0}$. After
this replacement, we may assume that

\begin{center}
$A=\overline{\Theta }_{\sigma }^{\ast }$.
\end{center}

\noindent By Lemma 3.2, $\overline{\Theta }_{\sigma }$ to a diffeomorphism $%
h_{\sigma }:X\rightarrow $ $X$. Therefore, replacing $f$ by

\begin{center}
$f^{\prime }=$ $f\circ h_{\sigma ^{-1}}:X\rightarrow X$
\end{center}

\noindent we obtain

\begin{center}
$A^{\prime }=\overline{\Theta }_{\sigma ^{-1}}^{\ast }\circ A=\mathrm{Id}$, 
\end{center}

\noindent and hence 

\begin{center}
$F^{\prime \ast }=\psi ^{\lambda }$. 
\end{center}

\noindent Since $\psi ^{\lambda }$ commutes with graded ring automorphisms,
undoing the replacement gives

\begin{center}
$f^{\ast }=\psi ^{\lambda }$ or $\psi ^{\lambda }\circ h_{\sigma }^{\ast }$, 
$\ \lambda =\deg (f)$,
\end{center}

\noindent where $\sigma \in \mathrm{Aut}(\Gamma _{G},\omega _{i})$. In
particular, the second possibility can occur only when 

\begin{center}
$\mathrm{Aut}(\Gamma _{G},\omega _{i})\neq \{\mathrm{id}\}$, 
\end{center}

\noindent and the corresponding Dynkin symmetry is the one classified in
Lemma 1.4. This completes the proof of Theorem B.\hfill $\square $

\section{Supplements to Theorem B}

For a generalized Grassmannian $X=G/P_{\{i\}}$ consider the degree-zero
component of the homotopy set $[X,X]$

\begin{center}
$[X,X]_{0}:=\left\{ [f]\in \lbrack X,X]\mid \deg (f)=0\right\} $.
\end{center}

\noindent Theorem B provides no information concerning the elements of $%
[X,X]_{0}$.

Letting $p:S(X)\rightarrow X$ denote the oriented circle bundle over $X$
whose Euler class is the K\"{a}hler class $\kappa _{X}\in H^{2}(X)$, we
obtain a sequence of fibrations

\begin{center}
$S^{1}\rightarrow S(X)\rightarrow X\overset{\kappa }{\rightarrow }\mathbb{C}%
P^{\infty }$,
\end{center}

\noindent where $\kappa $ denotes the classifying map of the K\"{a}hler
class $\kappa _{X}$. This sequence gives rise, for any topological space $Y$%
, to an exact sequence of homotopy sets

\begin{center}
$H^{1}(Y)=[Y,S^{1}]\rightarrow \lbrack Y,S(X)]\overset{p_{\ast }}{%
\rightarrow }[Y,X]\overset{\kappa _{\ast }}{\rightarrow }H^{2}(Y)$.
\end{center}

\noindent Taking $Y=X$, and using $H^{1}(X)=\{0\}$ together with

\begin{center}
$\kappa _{\ast }^{-1}(0)=[X,X]_{0}$,
\end{center}

\noindent we obtain the following result.

\bigskip

\noindent \textbf{Proposition 6.1.} For a generalized Grassmannian $X$, the
bundle map $p$ induces a bijection

\begin{center}
$p_{\ast }:[X,S(X)]\rightarrow \lbrack X,X]_{0}$, $[g]\rightarrow \lbrack
p\circ g]$.
\end{center}

This observation suggests that, in addition to the representation $R$ in
(1.1), one may introduce a secondary representation

\begin{center}
$R_{1}:[X,X]_{0}\rightarrow \mathrm{Hom}(H^{\ast }(S(X)),H^{\ast }(X))$,
\end{center}

\noindent which refines the original one. We point out that the cohomology
rings of both $X$ and the total space $S(X)$ can be effectively computed in
the framework of Schubert calculus (see, for example, \cite[Theorems 1-7:
Theorems 8-14]{DZ1}). Consequently, the set $\mathrm{Hom}(H^{\ast
}(S(X)),H^{\ast }(X))$ can, in principle, be explicitly determined.

Finally, regarding the geometric realization of the Adams operations $\psi
^{k}$ on the cohomology ring $H^{\ast }(X)$, as mentioned at the end of
Section 1, we refer to the following more general result of X. Z. Lin \cite[%
Theorem 1.3]{Lin}.

\bigskip

\noindent \textbf{Proposition 6.2. }Let $X=G/P$ be a flag manifold, and let $%
k$ be an integer coprime to the order of the Weyl group $W_{G}$. Then there
exists a self-map $f$ of $X$ such that

\begin{center}
$f^{\ast }=\psi ^{k}:H^{\ast }(X)\rightarrow H^{\ast }(X)$.
\end{center}

\noindent \textbf{Acknowledgement.} This work is supported by the National
Science Foundation of China under Grant nos. 12331003, 12288201; and by the
National Key R\TEXTsymbol{\backslash}\&D Program of China under Grant no.
2021YFA1002300.

\bigskip

Haibao Duan, dhb@math.ac.cn

Yau Mathematical Science Center, Tsinghua University, Beijing 100084;

Academy of Mathematics and Systems Sciences, Chinese Academy of Sciences,
Beijing 100190.

\bigskip

Ruizhi Huang, huangrz@amss.ac.cn

State Key Laboratory of Mathematical Sciences and Institute of Mathematics, Academy of Mathematics
and Systems Science, Chinese Academy of Sciences, Beijing 100190, China;

School of Mathematical Sciences, University of Chinese Academy of Sciences, Beijing 100049, China

\bigskip

Xuezhi Zhao, zhaoxve@mail.cnu.edu.cn

School of Mathematical Sciences, Capital Normal University, Beijing 100048.

\end{document}